\documentclass[12pt]{amsart}
\usepackage{hyperref}
\usepackage{booktabs}
\usepackage{amsmath}
\usepackage{amscd,amssymb}
\usepackage{epsfig}

\renewcommand{\proof}{\par\noindent{\it Proof.\ \ }}
\def\qed{\ifmmode\square\else\nolinebreak\hfill
$\Box$\fi\par\vskip12pt}

\newtheorem{thm}{Theorem}[section]
\newtheorem{prop}[thm]{Proposition}
\newtheorem{lem}[thm]{Lemma}

\newtheorem{cor}[thm]{Corollary}

\theoremstyle{definition}

\newtheorem{rem}[thm]{Remark}
\newtheorem{defn}[thm]{Definition}

\newcommand{\bC}{\mathbb C}
\newcommand{\bD}{\mathbb D}\newcommand{\bF}{\mathbb F}

\newcommand{\bQ}{\mathbb Q}

\newcommand{\bZ}{\mathbb Z}

\newcommand{\cM}{\mathcal M}

\def\ov{\overline}

\def\Aut{{\rm Aut}}

\def\Cay{{\rm Cay}}

 \def\PSL{{\rm PSL}}  
   \def\tr{{\rm tr}}
\def\Gal{{\rm Gal}}\def\Ker{{\rm Ker}}

\renewcommand{\>}{\rangle}

\begin{document}
\pagestyle{plain}

\begin{titlepage}

\begin{center}
\author{Yu Qing Chen and Tao Feng
}
\end{center}

\address{Department of Mathematics and Statistics, Wright State
University, Dayton, OH 45435}
\email{yuqing.chen@wright.edu}
\address{Department of Mathematics, Zhejiang University, Hangzhou 310027, China}
\email{tfeng@zju.edu.cn}

\begin{abstract}
In this paper, we present constructions of abelian Paley type group schemes by using multiplicative characters
of finite fields and Arasu-Dillon-Player difference sets. The constructions produce many new Paley type group schemes
that were previous unknown in our classification of Paley type group schemes in finite fields of small orders.

\end{abstract}

\title{Paley type group schemes from cyclotomic classes and Arasu-Dillon-Player difference sets}

\keywords{difference set; Paley group scheme; Paley type group scheme; Paley type partial difference set; skew Hadamard difference set; Singer difference set.
\\{\bf  Mathematics Subject Classification (2010) 05B10 · 05C25 · 05E18 · 05E30} }

\maketitle

\end{titlepage}

\section{Introduction}
A Paley type group scheme of a finite group $G$ is a 2-class association scheme $(G; R_0, R_1, R_2)$ obtained from a partition $D_0=\{1\}$, $D_1$, and $D_2$ of $G$ such that $D_1$ and $D_2$ satisfy the equation \begin{align}\label{def1}(1+2\bD^{(-1)})(1+2\bD)=|G|+(|G|-1)G\end{align} in the group ring $\bZ[G]$, where $\bD$ stands for the formal sum $$\sum_{x\in\bD}x$$ in the group ring $\bZ[G]$ of any subset $\bD$ in $G$ and the relations $R_i=\{(x,y)\in G\times G~|~xy^{-1}\in D_i\}$ for $i=0,1,2$. It is easy to show that $D_1$ satisfies (\ref{def1}) if and only if $D_2$ satisfies (\ref{def1}), and therefore we will also call any subset $\bD$ in $G$ that satisfies equation (\ref{def1}) a Paley type group scheme. Two Paley type group schemes $\bD_1$ and $\bD_2$ in the group $G$ are said to be equivalent if there is an automorphism $\alpha$ of $G$ such that $\bD_2=\bD_1^\alpha$. Associated with each Paley type group scheme $\bD$ in $G$ is a configuration $\mathfrak{C}(\bD)$ of $\bD$. When $|G|\equiv 1 \pmod 4$, a Paley type group scheme $\bD$ in $G$ is also known as a Paley type partial difference set. The configuration $\mathfrak{C}(\bD)$ is the Cayley graph $\Cay(G,\bD)$, which is a Paley type strongly regular graph and an example of Ramanujan graphs (see \cite{LPS} for definition). When $|G|\equiv 3 \pmod 4$, a Paley type group scheme $\bD$ in $G$ is also referred to as a skew Hadamard difference set. The configuration $\mathfrak{C}(\bD)$ is the Hadamard design ${\rm dev}(\bD)$ developed from $\bD$ (see \cite{BDL}). In this paper, we study Paley type group schemes in the additive group of finite fields.

Paley type group schemes were first studied by Paley~\cite{P}, who used quadratic residues in a finite field to construct
Hadamard matrices. The group schemes which are equivalent to quadratic residues in a finite field will be called Paley group schemes. Automorphism groups of configurations of Paley group schemes were determined by Carlitz~\cite{Car} and Kantor ~\cite{K0}. If we write $\Sigma_n$ for the symmetric group of degree $n$ and $\Gal(\bF_q)$ for the full automorphism group of the finite field $\bF_q$ of order $q$, then the results of Carlitz and Kantor concerning automorphism groups of configurations of Paley group schemes in \cite{Car, K0} can be summarized in the following theorem.
{\thm\cite[Theorem 8.1 and Corollary 8.2]{Car, K0} Let $q$ be a power of an odd prime, $\bF_q$ be a finite field of order $q$ and $S_{\bF_q}$ be the set of all non-zero quadratic residues in $\bF_q$. Then
$$\Aut(\mathfrak{C}(S_{\bF_q}))=\begin{cases}(\bF_q\rtimes S_{\bF_q})\rtimes\Gal(\bF_q)  &{\text if }\; q\notin\{3,7,11\},\\
                                             \PSL(2,q) &{\text if }\; q\in \{7,11\},\\
                                             \Sigma_3 &{\text if }\; q=3.\end{cases}$$}
If $q$ is a prime, then all Paley type group schemes in $\bF_q$ are Paley group schemes. If $q$ is a square of a prime, then all Paley type group schemes in $\bF_q$ can be constructed from partial congruence construction (see Theorem 2.2 in~\cite{M}). There are many construction methods of Paley type group schemes scattered in the literature. In~\cite{DA}, Davis discovered a product construction method and presented first family of Paley type group schemes in non-elementary abelian groups. Polhill~\cite{Pol1,Pol2} made tremendous advances in further developing product construction methods and constructed many Paley type group schemes in abelian groups which are not of prime power orders. Peisert~\cite[Theorem 3.1]{Pei} gave a classification of self-complementary symmetric graphs which produces a family of Paley type group schemes in $\bF_{p^{2n}}$ for $p\equiv 3\pmod 4$ by using the 4-class cyclotomic amorphic group scheme. Chen~\cite[Theorem 3.1]{Chen} also constructed a family of 4-class amorphic group schemes in $\bF_{q^4}$ which can be used to obtain Paley type group schemes. The discovery made by Ding and Yuan \cite{DY} re-energized research on Paley type group schemes in abelian groups of non-square orders. It is conjectured that such abelian groups must be elementary abelian $p$-groups and their exponent bounds were studied in~\cite{Ca,CXS,J,X}. Further new discoveries were made in \cite{DW} and Paley's construction were generalized in~\cite{WQ} by using ``quadratic residues" of commutative presemifields. Feng in~\cite{F} constructed a family of Paley type group schemes in extra-special $p$-groups of order $p^3$ and of exponent $p$ for $p>3$, and his construction was generalized by Chen and Polhill in~\cite{CP} using the flag group of finite fields (see \cite{C0}). In \cite{Mu}, Muzychuk obtained a large number of Paley type group schemes in $\bF_{q^3}$ and showed that the number of inequivalent such schemes grows exponentially. All results in \cite{CP,F,Mu} were generalized by Chen and Feng in \cite{CF}. In this paper we give a generalization of a cyclotomic construction of Paley type group schemes in \cite{FX} and present a new construction of Paley type group schemes in $\bF_{q^l}$ for odd $l$ by using Singer difference sets, Singer relative difference sets and Arasu-Dillon-Player difference sets (see Section~\ref{singer} for their definitions).

The constructions presented in this paper stem from our computer classification of Paley type group schemes in elementary abelian groups of small orders. All elementary abelian $p$-groups in this paper are presented as the additive group of finite fields and therefore all Paley type group schemes studied in this paper are in finite fields. By using MAGMA, we conducted exhaustive searches of all Paley type group schemes $\bD$ in finite fields of order $\leqslant240$ and classified their configurations up to isomorphism. We then compared the Paley type group schemes and their
configurations with those that can be constructed from known methods and our findings are tabulated in Table~\ref{t1}. The finite fields of order 243 and 343 are beyond the reach of exhaustive searches. But after we imposed some symmetric conditions on $\bD$, we managed to do complete searches for Paley type group schemes in these two fields that are invariant under the Galois group actions. The search results are exhibited in Table~\ref{t2}. The question marks in these tables indicate that there are Paley type group schemes found by computer which can not be constructed by all known methods prior to our constructions given in this paper.
It is interesting that the finite field of order 243 contains so many Galois invariant Paley type group schemes, and this may only be the tip of an iceberg because there may be many non-Galois invariant Paley type group schemes since they exist in other fields. For example, By Theorem 3.2 in~\cite{FX}, certain unions of cosets of the subgroup of order 95 in $\bF^*_{11^3}$ are Paley type group schemes. There are, besides Paley group scheme, 5 equivalence classes of such Paley type group schemes $\bD$, three of which have $|\Aut(\mathfrak{C}(\bD))|=3\cdot5\cdot11^3\cdot19$, and two of which have $|\Aut(\mathfrak{C}(\bD))|=5\cdot11^3\cdot19$. Clearly, the last two are not Galois invariant.

Throughout this paper $p$ is an odd prime, $q$ is a power of $p$, $l$ is a positive integer, $\bF$ is a finite field and $\bF_q$ is the finite field of order $q$, $S_{\bF_q}$ is the set of all non-zero quadratic residues of $\bF_q$, $N_{\bF_q}$ is the set of all quadratic non-residues of $\bF_q$, $\bF^*_q$ is the set of all  non-zero elements of $\bF_q$ which forms the multiplicative group of $\bF_q$, $\Gal(\bF_q)$ is the full automorphism group of the field $\bF_q$, and $\Gal(\bF_{q^l}/\bF_q)$ is the Galois group of $\bF_{q^l}$ over $\bF_q$. The main results of this paper are the following Theorems~\ref{main2}--\ref{main3}, and Theorem~\ref{main1} is a generalization of Theorem 3.2 in \cite{FX}.

Let $X$ be a subset of the quotient group $\bF^*_{q^l}/\bF^*_q$, and $\pi:\bF^*_{q^l}\to\bF^*_{q^l}/\bF^*_q$ be the natural projection homomorphism. We define the subset $\bD(X)$ in $\bF^*_{q^l}$ to be
\begin{align}\label{DX}\bD(X)=\{x\in S_{\bF_{q^l}}~|~\pi(x)\in X\}\cup\{x\in N_{\bF_{q^l}}~|~\pi(x)\notin X\}.\end{align}
Note that when $l$ is odd, the size $|\bD(X)|=(q^l-1)/2$ for all subsets $X$ in $\bF^*_{q^l}/\bF^*_q$.

{\thm\label{main2} Let $l$ be an odd integer and $X$ be a $((q^l-1)/(q-1),q^{l-1},q^{l-2}(q-1))$-difference set in $\bF^*_{q^l}/\bF^*_q$. Then $\bD(X)$ is a Paley type group scheme in $\bF_{q^l}$ if and only if $X$ is an Arasu-Dillon-Player difference set.}

The Paley type group schemes $\bD(X)$ in Theorem~\ref{main2} are $\Gal(\bF_{q^l})$ invariant when $p$ is a strong multiplier of the Arasu-Dillon-Player difference sets $X$.

{\thm\label{main1} Let $l$ be an odd integer, $n$ be a factor of $(q^l-1)/(q-1)$ such that $2$ is contained in the subgroup generated by $p$ in the multiplicative group $\bZ^*_n$ of the modular number ring $\bZ_n$, and $\gamma:\bF^*_{q^l}/\bF^*_q\to\bZ_n$ be the natural projection. Then for every subset $X$ in $\bZ_n$,  $\bD(\gamma^{-1}(X))$ is a Paley type group scheme in $\bF_{q^l}$.}

The Paley type group schemes in Theorem~\ref{main1} are unions of $n$ cyclotomic classes of order $2n$ in $\bF_{q^l}$ and Theorem 3.2 in \cite{FX} is an easy consequence of Theorem~\ref{main1}.
Theorem~\ref{main1} relaxed two very restrictive conditions of Theorem 3.2 in \cite{FX}, namely, the number $n$ in Theorem~\ref{main1} does not have to be a prime power and the index $[\bZ^*_n:\<p\>]$ in Theorem~\ref{main1} does not have to be $2$. For example, since the order of $2$ in $\bZ^*_{7\cdot31}$ is $15$, for every prime power $q\equiv 2\pmod {7\cdot31}$ and every odd integer $s$, there are Paley type group schemes in $\bF_{q^{15s}}$ that are unions of $7\cdot 31$ cyclotomic classes of order $2\cdot7\cdot31$. Theorem 3.6 and Corollary 3.7 in~\cite{FMX} and Theorem 3.6 in \cite{FX} will also be slightly generalized in Theorem~\ref{slight}.

By using an idea in the construction of Gordon-Mills-Welch difference sets as explained by Pott~\cite{Pott}, we obtain the following Theorem.

{\thm\label{main3} Let $s$ and $t$ be two positive odd integers and $\tilde{R}_{q^{st}/q^{t}}$ be the $((q^{st}-1)/(q^{t}-1), (q^t-1)/(q-1), q^{st-t},q^{st-2t}(q-1))$-Singer relative difference set in $\bF^*_{q^{st}}/\bF^*_q$ relative to $\bF^*_{q^t}/\bF^*_q$. For any subset $X$ in $\bF^*_{q^t}/\bF^*_q$, $\tilde{R}^{(-1)}_{q^{st}/q^{t}}X$ and $\tilde{R}^{(2)}_{q^{st}/q^{t}}X$ are subsets in $\bF^*_{q^{st}}/\bF^*_q$ and if $\bD(X)$ is a Paley type group scheme in $\bF_{q^{t}}$,
then $\bD(\tilde{R}^{(-1)}_{q^{st}/q^{t}}X)$ and $\bD(\tilde{R}^{(2)}_{q^{st}/q^{t}}X)$ are Paley type group schemes in $\bF_{q^{st}}$.}

The rest of the paper is organized as follows. In Section~\ref{singer}, we review difference sets, relative difference sets, weighing matrices, and Singer and Arasu-Dillon-Player difference sets. In Section~\ref{multi}, we give a multiplicative characterization of subsets $X$ in $\bF^*_{q^l}/\bF^*_q$ so that $\bD(X)$ are Paley type group schemes in $\bF_{q^l}$. Theorems~\ref{main2}--\ref{main3} are proved in Section~\ref{proof}. In the concluding section, Section~\ref{con}, we compare the Paley type group schemes constructed in this paper with those constructed from known methods in finite fields of small orders.

\begin{center}\begin{table}\caption
{Paley type group schemes $\mathbb{D}$ that are not Paley in finite fields $\bF$ of order $\leqslant 240$}\label{t1}
\begin{tabular}{|c|c|c|c|}
\toprule
$|\bF|$ & Number of inequivalent $\mathbb{D}$ and $\mathfrak{C}(\mathbb{D})$ & $|\Aut(\mathfrak{C}(\mathbb{D}))|$& Construction methods\\
\midrule
49&1&$2^3\cdot3^2\cdot7^2$&\cite[Theorem 3.4]{DA},\cite[Theorem 3.1]{Pei}\\
\midrule
81& 1& $2^5\cdot3^5\cdot5$&\cite[Theorem 3.1]{Chen},\cite[Theorem 3.1]{Pei} \\
\cline{2-4}
& $1$&$2^3\cdot3^5$&$?$\\
\cline{2-4}
& 1& $2^7\cdot3^4$&\cite[Theorem 3.4]{DA},\cite[Theorem 2.2]{M}\\
\midrule
121& 1& $2\cdot5^2\cdot11^2$&\cite[Theorem 3.4]{DA},\cite[Theorem 2.2]{M}\\
\cline{2-4}
& $1$&$2^2\cdot3\cdot5\cdot11^2$&\cite[Theorem 3.1]{Pei},\cite[Theorem 2.2]{M}\\
\cline{2-4}
& 1& $2^3\cdot5\cdot11^2$&\cite[Theorem 2.2]{M}\\
\midrule
125& 2 &$2^3\cdot3\cdot5^4$&\cite[Theorem 1.5]{CF}\\
\cline{2-4}
&1&$2^3\cdot3\cdot5^3$&?\\
\midrule
169&1&$2^3\cdot3^2\cdot13^2$&\cite[Theorem 3.4]{DA},\cite[Theorem 2.2]{M}\\
\cline{2-4}
&1&$2^2\cdot3^2\cdot13^2$&\cite[Theorem 2.2]{M}\\
\cline{2-4}
&2&$2^3\cdot3\cdot13^2$&\cite[Theorem 2.2]{M}\\
\bottomrule
\end{tabular}
\end{table}
\end{center}
\begin{center}
\begin{table}\caption{Non-Paley $\Gal(\bF)$ invariant Paley type group schemes $\mathbb{D}$ in finite fields $\bF$ of order 243 and 343}\label{t2}
\begin{tabular}{|c|c|c|c|}
\toprule
$|\bF|$ & Number of inequivalent $\mathbb{D}$ and $\mathfrak{C}(\mathbb{D})$ & $|\Aut(\mathfrak{C}(\mathbb{D}))|$& Construction methods\\
\midrule
243&1&$3^5\cdot5\cdot11$&\cite[Theorem 3.6]{FX}\\
\cline{2-4}
& 58& $3^5\cdot5$&?,\cite[Theorem 3.3]{DW},\cite[Corollary 3.7]{DY}\\
\midrule
343& $2$&$3^2\cdot7^4$&\cite[Theorem 1.5]{CF}\\
\cline{2-4}
& 2& $3^4\cdot7^3$&?\\
\cline{2-4}
& 1& $3^3\cdot7^3$&?\\
\cline{2-4}
& $7$&$3^2\cdot7^3$&$?$\\
\bottomrule
\end{tabular}
\end{table}
\end{center}

\section{Preliminaries on difference sets with Singer parameters}\label{singer}
In this section we review the basics of difference sets, relative difference sets, and
Singer weighing matrices. We then discuss difference sets with Singer parameters.
Some materials presented in this section and next section are taken from~\cite{C1}.

Given a finite group $G$ of order $v$, a $k$-subset $D$ of $G$ is
called a $(v,k,\lambda)$-difference set if for every $g\ne 1$ in
$G$, there are exactly $\lambda$ pairs of $(d_1, d_2)\in D\times D$
such that $d_1d_2^{-1}=g$. A $(v,k,\lambda)$-difference set in a
group $G$ gives rise to a $(v,k,\lambda)$-symmetric design whose
automorphism group contains $G$ as a subgroup which acts regularly
on both points and blocks of the design. If $G$ is a group of order $mn$ and
$N$ is a normal subgroup of $G$ of order $n$, a $k$-subset $D$ of $G$ is called
an $(m,n,k,\lambda)$-relative difference set relative to $N$ if for
every $g\in G\setminus N$, there are exactly $\lambda$ pairs
of $(d_1, d_2)\in D\times D$ such that $d_1d_2^{-1}=g$ and there is no
such expression for any $1\ne g\in N$. An
$(m,n,k,\lambda)$-relative difference set in a group $G$ relative to
a normal subgroup $N$ of $G$ gives rise to an
$(m,n,k,\lambda)$-symmetric divisible design whose automorphism
group contains $G$ as a subgroup which acts regularly on both points and blocks of the design and the normal subgroup $N$ is the
stabilizer in $G$ of the point classes and parallel classes of the
design.  when $n=1$, an $(m,n,k,\lambda)$-relative difference set
is simply an $(m,k,\lambda)$-difference set. If $D$ is an
$(m,n,k,\lambda)$-relative difference set in a group $G$ relative to
a normal subgroup $N$ of $G$, an automorphism $\sigma\in\Aut(G)$
such that $\{\sigma(g)~|~g\in D\}=\{gh~|~g\in D\}$ for some
$h\in G$ is called a multiplier of $D$. Multipliers of $D$ form a
group, which will be called the multiplier group of $D$ and will be denoted
by $\mathcal{M}(D)$. Each multiplier in $\mathcal{M}(D)$ induces an
automorphism of the design obtained from $D$. We call the subgroup $\cM_0(D)=\{\sigma\in\Aut(G)~|~\sigma(g)\in D\text{ for all }g\in D\}$
of $\cM(D)$ the strong multiplier group of $D$. For more details on
difference sets and relative difference sets, we refer the reader to Beth et al.
\cite{BDL} and Pott \cite{Pott}.

Let $G$ be a finite group. The group ring $\bZ[G]$ of the group
$G$ is the set of formal sums
$$\sum_{g\in G}a_gg,$$ where $a_g\in\bZ$ is the integer coefficient of $g$ in
the formal sum, endowed with the addition $$\left(\sum_{g\in
G}a_gg\right)+\left(\sum_{g\in G}b_gg\right)=\sum_{g\in G}\left(a_g+b_g\right)g,$$ and
multiplication $$\left(\sum_{g\in G}a_gg\right)\left(\sum_{g\in G}b_gg\right)=\sum_{g\in
G}\left(\sum_{h\in G}a_{gh^{-1}}b_{h}\right)g.$$ It is clear that $\bZ[G]$ is a
ring. For any subset $X$ of $G$ we often identify $X$ with
the group ring element
$$X=\sum_{g\in X}g$$ and for any group ring element $$A=\sum_{g\in G}a_gg\in\bZ[G],$$
we define $A^{(t)}$ to be
$$A^{(t)}=\sum_{g\in G}a_gg^t\in\bZ[G]$$ for any integer $t\in\bZ$. Using the group ring
notation, a $k$-element subset $D$ of $G$ is a $(v,k,\lambda)$-
difference set if and only if
$$DD^{(-1)}=(k-\lambda)+\lambda G$$ in $\bZ[G]$, or an $(m,n,k,\lambda)$-relative
difference set relative to a normal subgroup $N$ of $G$ if
and only if
$$DD^{(-1)}=k+\lambda (G-N)$$ in $\bZ[G]$.
Let $R$ be a ring and $R^*$ be the multiplicative group of the
invertible elements of $R$. If $f:G\rightarrow R^*$ is a group
homomorphism, then $f$ induces a ring homomorphism
$f:\bZ[G]\rightarrow R$ by $\bZ$-linearly extending $f$ from $G$ to
$\bZ[G]$. Therefore for any abelian group $G$, every character
$\chi:G\to\bC^*$ of $G$ can be extended to a ring homomorphism
$\chi:\bZ[G]\to\bC$, where $\bC$ is the complex number field.  We
denote the set of all characters of $G$ by $\widehat{G}$ and
$\widehat{G}$ forms an abelian group under point-wise multiplication as
functions from $G$ to $\bC^*$. The trivial homomorphism from $G$ to
$\bC^*$ is called the principal character of $G$. We call $\widehat{G}$
the dual of $G$. By the Fourier inversion formula, one has

{\lem\label{Dset} A $k$-element subset $D$ of an abelian group $G$
of order $v$ is a $(v,k,\lambda)$-difference set if and only if for
every non-principal character $\chi\in\widehat{G}$,
$|\chi(D)|=\sqrt{k-\lambda}$.

A $k$-element subset $D$ of an abelian group $G$ of order $mn$ is an
$(m,n,k,\lambda)$-relative difference set relative to a subgroup $N$
of $G$ of order $n$ if and only if for every character $\chi\in\widehat{G}$
which is non-trivial on $N$, $|\chi(D)|=\sqrt{k}$, and for every
non-principal character $\chi\in\widehat{G}$ which is trivial on $N$,
$|\chi(D)|=\sqrt{k-\lambda n}$.}\vspace{2mm}

If $G$ is an elementary abelian $p$-group, we can identify $G$ with
the additive group of a finite field
$\bF$ of the same size. Let $\xi_p$ be a primitive $p$-th root of unity and
$\tr:\bF\to\bF_p$ be the trace map, that is $$\tr(x)=\sum_{\sigma\in\Gal(\bF)}x^\sigma$$ for all $x\in\bF$. Then
every character in $\widehat{\bF}$ is given by
\begin{eqnarray*}&\chi_\alpha:\bF\to\bC^*&\\
&\chi_\alpha(x)=\xi_p^{\tr(\alpha x)}& \mbox{ for all }x\in\bF,
\end{eqnarray*} where $\alpha\in\bF$.

In order to obtain a multiplicative description of Paley type group schemes in
additive groups of finite fields, we need the Singer difference sets and Singer
relative difference sets.  Let
$\tr_{q^l/q}:\bF_{q^l}\rightarrow\bF_q$ be the relative trace
map from $\bF_{q^l}$ to $\bF_q$, i.e.
\begin{align*}
\tr_{q^l/q}(x)&=\sum_{\sigma\in\Gal(\bF_{q^l}/\bF_q)}x^\sigma
\end{align*} for all
$x\in\bF_{q^l}$. Let
$R_{q^l/q}=\{x\in\bF_{q^l}^{*}~|~\tr_{q^l/q}(x)=1\}\subset\bF_{q^l}^{*}$ and
$S_{q^l/q}\subset \bF_{q^l}^{*}/\bF_q^{*}$ be the image of $R_{q^l/q}$ under the
natural projection map from $\bF_{q^l}^{*}$ to
$\bF_{q^l}^{*}/\bF_q^{*}$. The set $R_{q^l/q}$ is the
$((q^l-1)/(q-1),q-1,q^{l-1},q^{l-2})$-Singer relative difference set
in $\bF_{q^l}^{*}$ relative to $\bF_q^*$. The set $S_{q^l/q}$ is the
$((q^l-1)/(q-1),q^{l-1},q^{l-2}(q-1))$-Singer difference set in
$\bF_{q^l}^{*}/\bF_q^{*}$. The complement of $S_{q^l/q}$ in
$\bF_{q^l}^{*}/\bF_q^{*}$ is the
$((q^l-1)/(q-1),(q^{l-1}-1)/(q-1),(q^{l-2}-1)/(q-1))$-Singer
difference set and it can be obtained from the projection of the
hyperplane $\{x\in\bF_{q^l}^{*}~|~\tr_{q^l/q}(x)=0\}$ in
$\bF_{q^l}^{*}$ to the quotient group $\bF_{q^l}^{*}/\bF_q^{*}$.
The prime $p$ is in the strong multiplier groups $\cM_0(R_{q^l/q})$ and $\cM_0(S_{q^l/q})$.
When $l=st$ for some positive integers $s$ and $t$, $R_{q^{st}/q}=R_{q^{st}/q^t}R_{q^t/q}$
as $\tr_{q^{st}/q}=\tr_{q^t/q}\circ\tr_{q^{st}/q^t}$, and the Singer difference set
$S_{q^{st}/q}=\tilde{R}_{q^{st}/q^t}S_{q^t/q}$, where $\tilde{R}_{q^{st}/q^t}$ is the image
of $R_{q^{st}/q^t}$ under the natural projection $\bF^*_{q^{st}}\to\bF^*_{q^{st}}/\bF^*_q$.
We call this decomposition of $S_{q^{st}/q}$ the Gordon-Mills-Welch decomposition which forms the
foundation of the Gordon-Mills-Welch construction~\cite{GMW} of difference sets having same parameters
as that of $S_{q^{st}/q}$. In~\cite{Pott}, Pott presented a
construction which is more general than the Gordon-Mills-Welch construction.

{\prop\label{GGMW}\cite[Proposition 3.2.1]{Pott} Let $R$ be an $(m,n,k,\lambda)$-relative difference
set in an abelian group $G$ relative to a subgroup $N$, and $T$ be an $(n/n',n',k',\lambda')$-relative
difference set in $N$ relative to a subgroup $N'$ of $N$. If $k\lambda'-k'\lambda=\lambda\lambda'(n-n')$,
then the subset $RT$ is an $(mn/n',n',kk',k\lambda')$-relative difference set in $G$ relative to $N'$. }

In \cite{AC}, the following product formula for difference sets with $k-\lambda$ dividing $\lambda$ was proved.

{\thm\label{product}\cite[Theorem 2.3]{AC} Let $G$ be a group of order $v$ and $D_1$,
$D_2$, $\cdots$, $D_{2r+1}$ be $(v,k,\lambda)$-difference sets in
$G$ with $n|\lambda$, where $n=k-\lambda$. If $n^r$ divides $D_1D_2\cdots D_{2r+1}$ in
$\bZ[G]$, then there is a $(v,k,\lambda)$ difference set $D$ in $G$
such that
\begin{align}\label{pdt}D_1D_2\cdots
D_{2r+1}=n^r\left(\frac{k[(1+sv)^r-1]}{v}G+D\right)=(n+\lambda
G)^rD\end{align} in $\bZ[G]$, where $s=\lambda/n$.}

Theorem~\ref{product} can also be formulated as

{\thm\label{product1} Let $G$ be a group of order $v$ and $D_1$,
$D_2$, $\cdots$, $D_{2r+1}$ be $(v,k,\lambda)$ difference sets in
$G$ with $n|\lambda$, where $n=k-\lambda$. If $D_1D_2\cdots D_r$ divides $D_{r+1}D_{r+2}\cdots D_{2r+1}$ in
$\bZ[G]$, then there is a $(v,k,\lambda)$ difference set $D$ in $G$
such that
\begin{align}\label{pdt1}D_1D_2\cdots D_rD=D_{r+1}D_{r+2}\cdots D_{2r+1}\end{align} in $\bZ[G]$.}

The product formula given in equation~(\ref{pdt1}) resembles in some way the non-unique factorization of algebraic integers in number fields. For instance, all $((q^l-1)/(q-1),q^{l-1},q^{l-2}(q-1))$-difference sets
satisfy the condition that $k-\lambda$ divides $\lambda$, and Arasu, Dillon and Player in~\cite{AD} obtained many different factorizations of $S_{q^l/q}^{(-1)}S_{q^l/q}^{(2)}$ in $\bZ[\bF^*_{q^l}/\bF^*_q]$. These factorizations
promote the following definition.

{\defn\label{adpd} A $((q^l-1)/(q-1),q^{l-1},q^{l-2}(q-1))$-difference set $A$ in $\bF^*_{q^l}/\bF^*_q$ will be called an Arasu-Dillon-Player difference set if $A$ divides $S_{q^l/q}^{(-1)}S_{q^l/q}^{(2)}$ in $\bZ[\bF^*_{q^l}/\bF^*_q]$.
{\rem\label{remadp}  There are actually several different but equivalent definitions of Arasu-Dillon-Player difference sets. For example, a $((q^l-1)/(q-1),q^{l-1},q^{l-2}(q-1))$-difference set $A$ in $\bF^*_{q^l}/\bF^*_q$ is an Arasu-Dillon-Player difference set if and only if $S_{q^l/q}^{(2)}$ divides $AS_{q^l/q}$ in $\bZ[\bF^*_{q^l}/\bF^*_q]$.}

By Theorem~\ref{product1}, Arasu-Dillon-Player difference sets are in pairs, that is if $A$ is an Arasu-Dillon-Player difference set in $\bF^*_{q^l}/\bF^*_q$, then there is another Arasu-Dillon-Player difference set $B$ in $\bF^*_{q^l}/\bF^*_q$ such that $AB=S_{q^l/q}^{(-1)}S_{q^l/q}^{(2)}$ in $\bZ[\bF^*_{q^l}/\bF^*_q]$. We will call the pair $A$ and $B$ a dual pair of Arasu-Dillon-Player difference sets, and call $A$ the dual of $B$. The following theorem can be found in \cite{A0} and \cite{AD}.

{\thm\label{adp}~\cite[Theorem 6.17 and 6.19]{A0} If $q^r+1$ is prime to $(q^l-1)/(q-1)$, then $S^{(1+q^r)}_{q^l/q}$ is an Arasu-Dillon-Player difference set. If $r$ is prime to $l$, then $S^{(\frac{1+3^r}{2})}_{3^l/3}$ is an Arasu-Dillon-Player difference set.}

We write $A_{q^l/q}(1+q^r)$ and $A_{3^l/3}(\frac{1+3^r}{2})$ for the dual Arasu-Dillon-Player difference sets of $S^{(1+q^r)}_{q^l/q}$  and $S^{(\frac{1+3^r}{2})}_{3^l/3}$ respectively. These difference set will be used in the last section to obtain new Paley type group schemes.

Using the Gordon-Mills-Welch decomposition $S_{q^{st}/q}=\tilde{R}_{q^{st}/q^t}S_{q^t/q}$ and Proposition~\ref{GGMW}, one can construct more Arasu-Dillion-Player difference sets.

{\thm If $A$ is a $((q^t-1)/(q-1),q^{t-1},q^{t-2}(q-1))$-Arasu-Dillion-Player difference set in $\bF^*_{q^t}/\bF^*_q$, then $\tilde{R}^{(-1)}_{q^{st}/q^t}A$ and $\tilde{R}^{(2)}_{q^{st}/q^t}A$ are $((q^{st}-1)/(q-1),q^{st-1},q^{st-2}(q-1))$-Arasu-Dillion-Player difference sets in $\bF^*_{q^{st}}/\bF^*_q$.}
\proof It is easy to check that the parameters of $\tilde{R}_{q^{st}/q^t}$ and $A$ satisfy the condition in Proposition~\ref{GGMW} with $n'=1$. Since $A$ is a $((q^t-1)/(q-1),q^{t-1},q^{t-2}(q-1))$-Arasu-Dillion-Player difference set in $\bF^*_{q^t}/\bF^*_q$, the element $A$ divides $S_{q^t/q}^{(-1)}S_{q^t/q}^{(2)}$ in $\bZ[\bF^*_{q^t}/\bF^*_q]\subseteq\bZ[\bF^*_{q^{st}}/\bF^*_q]$ and therefore both $\tilde{R}^{(-1)}_{q^{st}/q^t}A$ and $\tilde{R}^{(2)}_{q^{st}/q^t}A$ divide $\tilde{R}^{(-1)}_{q^{st}/q^t}\tilde{R}^{(2)}_{q^{st}/q^t}S_{q^t/q}^{(-1)}S_{q^t/q}^{(2)}=S_{q^{st}/q}^{(-1)}S_{q^{st}/q}^{(2)}$ in $\bZ[\bF^*_{q^{st}}/\bF^*_q]$.\qed

Another object that has relevance to our construction is the Singer
circulant weighing matrix. Given two positive integers $k$ and $n$,
an $n\times n$ matrix $M=(m_{i,j})_{n\times n}$ is called a $(k,n)$
weighing matrix if $m_{i,j}^3=m_{i,j}$ for all $i,j=1,2,\cdots,n$
and $MM^\top=kI_n$, where $M^\top$ is the transpose of $M$ and $I_n$
is the $n\times n$ identity matrix. The set of all $(k,n)$ weighing
matrices is denoted by ${\rm W}(k,n)$. An $n\times n$ matrix
$M=(m_{i,j})_{n\times n}$ is said to be circulant if
$m_{i',j'}=m_{i,j}$ whenever $j'-i'\equiv j-i$ $(\!\!\!\mod n)$. The set
of all circulant  $(k,n)$ weighing matrices is denoted by ${\rm CW}(k,n)$.
Circulant matrices with integer entries can be viewed as group ring
elements of a cyclic group. This is because if $M=(m_{i,j})_{n\times
n}$ is a circulant matrix, let $a_i=m_{1,i+1}$ for
$i=0,1,\cdots,n-1$, $g=(g_{i,j})_{n\times n}$ with
\[g_{i,j}=\left\{\begin{array}{ll}1,&\mbox{if } j-i\equiv 1~(\!\!\!\!\!\mod n)\\
0&\mbox{otherwise, }\end{array}\right.\] and $G=\<g\>\cong\bZ_n$,
then $M=a_0g^0+a_1g^1+a_2g^2+\cdots+a_{n-1}g^{n-1}\in\bZ[G]$ and
$M^\top=a_0g^0+a_1g^{-1}+a_2g^{-2}+\cdots+a_{n-1}g^{-(n-1)}=M^{(-1)}\in\bZ[G]$.
Hence a matrix $M\in{\rm CW}(k,n)$ simply means that $M$ is an element in
$\bZ[G]$ with $-1$, $0$, $1$ coefficients and $MM^{(-1)}=k$ in
$\bZ[G]$. Circulant weighing matrices were studied extensively in
\cite{A, AL, AL2, AM, LMS}. If $q$ is a power of a prime and
$l$ is an odd positive integer, then the image $R$ in
$\bF^*_{q^l}/S_{\bF_q}$ of the Singer relative difference set $R_{q^l/q}$ in
$\bF^*_{q^l}$ relative to $\bF^*_q$ via the natural projection
$\bF^*_{q^l}\to \bF^*_{q^l}/S_{\bF_q}$ is the
$((q^l-1)/(q-1),2,q^{l-1},q^{l-2}(q-1)/2)$-Singer relative difference set
in $\bF_{q^l}^{*}/S_{\bF_q}$ relative to $\bF_q^*/S_{\bF_q}$.  The
group $\bF_{q^l}^{*}/S_{\bF_q}\cong
(\bF_{q^l}^{*}/\bF^{*}_q)\times(\bF_q^*/S_{\bF_q})$. If we replace the
non-trivial element of $\bF_q^*/S_{\bF_q}\cong\bZ_2$ with $-1$, then the relative
difference set
$R$ becomes a $(q^{l-1},(q^l-1)/(q-1))$ circulant
weighing matrix. We call this matrix the Singer circulant weighing
matrix and will denote it by $W_{q^l/q}$.

\section{A multiplicative characterization of certain Paley type group schemes
in finite fields}\label{multi}
We now discuss a characterization of certain Paley type group schemes in finite fields by
using the multiplicative group of the fields. The materials presented here are contained in
\cite{C1} for $q\equiv 3\pmod 4$. The idea was used by Dillon in \cite{D} and \cite{D1} for Hadamard difference sets in
elementary abelian $2$-groups, which is the reverse of the method used by Dillon in \cite{D2}.

{\defn A subset $X$ in $\bF^*_{q^l}$ will be called a {\it projective half-point set} over $\bF_q$ if
$X$ is invariant under the multiplication of $S_{\bF_q}$, i.e. $X$ is a union of cosets of
$S_{\bF_q}$, and the intersection of $X$ with every coset of $\bF^*_q$ is of size $(q-1)/2$.}

If we view $\bF^*_{q^l}$ as the projective space ${\rm PG}(l-1,\bF_q)$, then a coset of $\bF^*_q$
is a point in ${\rm PG}(l-1,\bF_q)$ and the set $X$ consists of half of each point in
${\rm PG}(l-1,\bF_q)$. The following theorem gives a necessary and sufficient condition for a
projective half-point set in $\bF^*_{q^l}$ over $\bF_q$ to be a Paley type group scheme in $\bF_{q^l}$.

{\thm\label{mul1} Let $D$ be a projective half-point set in $\bF^*_{q^l}$ over $\bF_q$. The subset
$D$ of $\bF^*_{q^l}$ is a Paley type group scheme in the additive group of $\bF_{q^l}$ if and only
if there is a subset $\widehat{D}$ in $\bF^*_{q^l}$ such that
$D$ and $\widehat{D}$ satisfy the equation
\begin{align}\label{mul11}D^{(-1)}R_{q^l/q}=q^{(l-1)/2}\widehat{D}
+\frac{q^{l-1}-q^{(l-1)/2}}{2}\bF^*_{q^l}\end{align} in the group ring
$\bZ[\bF^*_{q^l}]$, where $R_{q^l/q}$ is the
$((q^l-1)/(q-1),q-1,q^{l-1},q^{l-2})$-Singer relative difference set
in $\bF^*_{q^l}$ relative to $\bF^*_q$.}
\proof By the Fourier inversion formula, we only need to show that equation (\ref{mul11}) is equivalent to
$$|\chi_g(1+2D)|=\sqrt{q^l}$$ for all $g\in\bF^*_{q^l}$. The group $\bF_{q^l}^*$ can be partitioned into three
subsets \begin{align*}
H^*&=\{x\in\bF_{q^l}^*~|\tr_{q^l/q}(x)=0\},\\
S_{\bF_q}R_{q^l/q}&=\{x\in\bF_{q^l}^*~|\tr_{q^l/q}(x)\in S_{\bF_q}\},\\
N_{\bF_q}R_{q^l/q}&=\{x\in\bF_{q^l}^*~|\tr_{q^l/q}(x)\in N_{\bF_q}\}.
\end{align*} Given an element $g\in\bF_{q^l}^*$, the character sum
\begin{align*}\chi_g(1+2D)=&1+2\sum_{x\in D}\xi_p^{\tr_{q^l/p}(gx)}\\=&
1+2\sum_{x\in D,gx\in H^*}\xi_p^{\tr_{q^l/p}(gx)}
+2\sum_{x\in D,gx\in S_{\bF_q}R_{q^l/q}}\xi_p^{\tr_{q^l/p}(gx)}+2\sum_{x\in
D,gx\in N_{\bF_q}R_{q^l/q}}\xi_p^{\tr_{q^l/p}(gx)}
\end{align*}
with \begin{align*} 1+2\sum_{x\in D,gx\in H^*}\xi_p^{\tr_{q^l/p}(gx)}
=&1+2\sum_{x\in D,gx\in H^*}\xi_p^{\tr_{q/p}(\tr_{q^l/q}(gx))}\\
=&1+2|gD\cap H^*|=1+|H^*|=q^{l-1},
\end{align*}
\begin{align*}
2\sum_{x\in D,gx\in S_{\bF_q}R_{q^l/q}}\xi_p^{\tr_{q^l/p}(gx)}
=&2\sum_{x\in D,gx\in S_{\bF_q}R_{q^l/q}}\xi_p^{\tr_{q/p}(\tr_{q^l/q}(gx))}\\
=&2|gD\cap R_{q^l/q}|\sum_{x\in S_{\bF_q}}\xi_p^{\tr_{q/p}(x)}\\=&2|gD\cap R_{q^l/q}|\chi_1(S_{\bF_q}),
\end{align*} and
\begin{align*}
2\sum_{x\in D,gx\in N_{\bF_q}R_{q^l/q}}\xi_p^{\tr_{q^l/p}(gx)}
=&2\sum_{x\in D,gx\in N_{\bF_q}R_{q^l/q}}\xi_p^{\tr_{q/p}(\tr_{q^l/q}(gx))}\\
=&2|gD\cap\omega R_{q^l/q}|\sum_{x\in N_{\bF_q}}\xi_p^{\tr_{q/p}(x)}\\
=&2|gD\cap\omega R_{q^l/q}|\chi_1(N_{\bF_q}),\\
\end{align*} where $\omega$ is a primitive element in $\bF_q$.  Note that $|gD\cap\omega R_{q^l/q}|=|g\omega^{-1}D\cap R_{q^l/q}|$.
Since $D\cap\omega^{-1}D=\emptyset$ and $D\cup\omega^{-1}D=\bF^*_{q^l}$ as $D$ is a projective half-point set in $\bF^*_{q^l}$ over $\bF_q$,
we also have $gD\cap g\omega^{-1}D=\emptyset$ and $gD\cup g\omega^{-1}D=\bF^*_{q^l}$ for all $g\in\bF^*_{q^l}$. This implies that
\begin{align}\label{plus}|gD\cap R_{q^l/q}|+|gD\cap\omega R_{q^l/q}|=|gD\cap R_{q^l/q}|+|g\omega^{-1}D\cap R_{q^l/q}|=|R_{q^l/q}|=q^{l-1}.\end{align}
Therefore $D$ is a Paley type group scheme if and only if
\begin{align*}
\sqrt{q^l}=|\chi_g(1+2D)|=&||gD\cap R_{q^l/q}|\chi_1(1+2S_{\bF_q})+|gD\cap\omega R_{q^l/q}|\chi_1(1+2N_{\bF_q})|\\
=&|(|gD\cap R_{q^l/q}|-|gD\cap\omega R_{q^l/q}|)\chi_1(1+2S_{\bF_q})|,\end{align*}
which is equivalent to
\begin{equation}\label{minus}|gD\cap R_{q^l/q}|-|gD\cap\omega R_{q^l/q}|=\pm q^{\frac{l-1}{2}}.\end{equation}
Combining (\ref{plus}) and (\ref{minus}), we get
\begin{align}|gD\cap R_{q^l/q}|=&\frac{q^{l-1}\pm q^{\frac{l-1}{2}}}{2},\label{+}\end{align}
Let
$$\widehat{D}=\{g\in\bF^*_{q^l}~|~|gD\cap R_{q^l/q}|=\frac{q^{l-1}+
q^{\frac{l-1}{2}}}{2}\}=\{g\in\bF^*_{q^l}~|~\chi_g(1+2D)=q^{\frac{l-1}{2}}\chi_1(1+2S_{\bF_q})\}.$$
Then
$$\omega\widehat{D}=\{g\in\bF^*_{q^l}~|~|gD\cap R_{q^l/q}|=\frac{q^{l-1}-
q^{\frac{l-1}{2}}}{2}\}=\{g\in\bF^*_{q^l}~|~\chi_g(1+2D)=q^{\frac{l-1}{2}}\chi_1(1+2N_{\bF_q})\}$$
and
\begin{align*}D^{(-1)}R_{q^l/q}&=\frac{q^{l-1}+q^{\frac{l-1}{2}}}{2}\widehat{D}+\frac{q^{l-1}-
q^{\frac{l-1}{2}}}{2}(\omega\widehat{D})\\&=q^{\frac{l-1}{2}}\widehat{D}+\frac{q^{l-1}-q^{\frac{l-1}{2}}}{2}\bF^*_{q^l}\end{align*} in the group
ring $\bZ[\bF^*_{q^l}]$. \qed

{\rem\label{mul1r} From the proof it is easy to see that $\widehat{D}$ is also a projective-half point set and satisfies
$$\widehat{D}^{(-1)}R_{q^l/q}=q^{\frac{l-1}{2}}D
+\frac{q^{l-1}-q^{\frac{l-1}{2}}}{2}\bF^*_{q^l}$$ in the group
ring $\bZ[\bF^*_{q^l}]$. Therefore $\widehat{D}^{(-1)}$ is also a Paley type group scheme and it is the dual of $D$.}

The next theorem is equivalent to Theorem~\ref{mul1} but much easier to use. Its proof requires the following lemma of Ma~\cite{Ma}.

{\lem\cite[Lamma 3.4]{Ma}\label{sqsum} Let $a_1$, $a_2$, $\cdots$, $a_m$ be integers
and
$$\sum_{i=1}^{m}a_i=n.$$ If $n=qm+r$, where $q$ and $r$ are integers
and $0\leqslant r<m$, then $$\sum_{i=1}^{m}a^2_i\geqslant
(m-r)q^2+r(q+1)^2$$ and equality holds if and only if
$|a_i-a_j|\leqslant 1$ for all $1\leqslant i,j\leqslant m,$ i.e.
there are exactly $m-r$ of $a_1$, $a_2$, $\cdots$, $a_m$ with value
$q$ and the remaining $r$ of them with value $q+1$. }\vspace{1.0mm}

\proof  If there are $i$ and $j$ such that $a_i-a_j>1$, let
$a{'}_k=a_k$ for $k\ne i$ or $j$, $a{'}_i=a_i-1$ and
$a{'}_j=a_j+1$, then
$$\sum_{k=1}^{m}a{'}_k=\sum_{k=1}^{m}a_k=n$$ and
$$\sum_{k=1}^{m}{a'}^2_k=\sum_{k=1}^{m}{a}^2_k+2(1+a_j-a_i)<\sum_{k=1}^{m}{a}^2_k.$$
Hence $\sum_{k=1}^{m}{a}^2_k$ attains minimum if and only if $m-r$
of the integers $a_1$, $a_2$, $\cdots$, $a_m$ are equal to $q$ and
the remaining $r$ of them are equal to $q+1$.\qed

{\thm\label{mul2} A projective half-point set $D$ in $\bF^*_{q^l}$ over $\bF_q$ is a Paley type group scheme in $\bF_{q^l}$ if and only if $D^{(-1)}R_{q^l/q}$ is divisible by $q^{(l-1)/2}$ in
the group ring $\bZ[\bF^*_{q^l}]$, where $R_{q^l/q}$ is the
$((q^l-1)/(q-1),q-1,q^{l-1},q^{l-2})$-Singer relative difference set
in $\bF^*_{q^l}$ relative to $\bF^*_q$.}
\proof If $D$ is a Paley type group scheme in $\bF_{q^l}$, then by Theorem~\ref{mul1}, $D^{(-1)}R_{q^l/q}$ is divisible by $q^{(l-1)/2}$ in the group ring $\bZ[\bF^*_{q^l}]$. Conversely, if
$D^{(-1)}R_{q^l/q}$ is divisible by $q^{(l-1)/2}$ in the group ring $\bZ[\bF^*_{q^l}]$, then
$$\frac{D^{(-1)}R_{q^l/q}}{q^{\frac{l-1}{2}}}=\sum_{g\in\bF^*_{q^l}}a_gg\in\bZ[\bF^*_{q^l}]$$ for some integers $a_g$ and
\begin{align}\sum_{g\in\bF^*_{q^l}}a_g=\frac{q^{\frac{l-1}{2}}(q^l-1)}{2}=\frac{q^{\frac{l-1}{2}}-1}{2}(q^l-1)+\frac{q^l-1}{2}.\end{align} Since
\begin{align*}D^{(-1)}R_{q^l/q}(D^{(-1)}R_{q^l/q})^{(-1)}&=DD^{(-1)}(q^{l-1}+q^{l-2}(\bF^*_{q^l}-\bF^*_{q}))\\
                                                         &=q^{l-1}DD^{(-1)}+\frac{q^{l-2}(q^l-1)^2}{4}\bF^*_{q^l}-\frac{q^{l-2}(q^l-1)(q-1)}{4}\bF^*_{q^l}\\
                                                         &=q^{l-1}DD^{(-1)}+\frac{q^{l-1}(q^l-1)(q^{l-1}-1)}{4}\bF^*_{q^l}\end{align*}
as $D$ is a projective half-point set in $\bF^*_{q^l}$ over $\bF_q$, we find that
\begin{align}\sum_{g\in\bF^*_{q^l}}a_g^2=\frac{q^l-1}{2}+\frac{(q^l-1)(q^{l-1}-1)}{4}=\left(\frac{q^{\frac{l-1}{2}}-1}{2}\right)^2\frac{q^l-1}{2}+\left(\frac{q^{\frac{l-1}{2}}+1}{2}\right)^2\frac{q^l-1}{2}.\end{align}
By Lemma~\ref{sqsum}, there is a subset $\widehat{D}$ of size $(q^l-1)/2$ in $\bF^*_{q^l}$ such that
$$\frac{D^{(-1)}R_{q^l/q}}{q^{\frac{l-1}{2}}}=\widehat{D}+\frac{q^{\frac{l-1}{2}}-1}{2}\bF^*_{q^l}.$$ By Theorem~\ref{mul1}, $D$ is a Paley type group scheme.\qed

Let $\mu:\bF^*_{q^l}\to\bF^*_{q^l}/S_{\bF_q}$ be the natural projection map. There is an one-to-one correspondence between projective half-point sets in $\bF^*_{q^l}$ over $\bF^*_{q}$ and transversals of $\bF^*_{q}/S_{\bF_q}$ in $\bF^*_{q^l}/S_{\bF_q}$. Combining Theorems~\ref{mul1}, Remark~\ref{mul1r} and Theorem~\ref{mul2}, one has

{\cor Let $D$ be a transversal of the subgroup $\bF^*_q/S_{\bF_q}$
in $\bF^*_{q^l}/S_{\bF_q}$ and $R$ be the
$((q^l-1)/(q-1),2,q^{l-1},q^{l-2}(q-1)/2)$-Singer relative
difference set in $\bF^*_{q^l}/S_{\bF_q}$ relative to
$\bF^*_q/S_{\bF_q}$. Then the following statements are equivalent:
\begin{itemize}
\item[(i)] $\mu^{-1}(D)$ is a Paley type group scheme in the
additive group of $\bF_{q^l}$;
\item[(ii)] there is another transversal $\widehat{D}$ of
the subgroup $\bF^*_q/S_{\bF_q}$
in $\bF^*_{q^l}/S_{\bF_q}$ such that
$$D^{(-1)}R=q^{(l-1)/2}\widehat{D}
+\frac{q^{l-1}-q^{(l-1)/2}}{2}\bF^*_{q^l}/S_{\bF_q}$$ in
$\bZ[\bF^*_{q^l}/S_{\bF_q}]$;
\item[(iii)] $D^{(-1)}R$ is divisible by
$q^{(l-1)/2}$ in $\bZ[\bF^*_{q^l}/S_{\bF_q}]$.\end{itemize}}\vspace{2mm}

When $l$ is odd, the group $\bF^*_{q^l}/S_{\bF_q}\cong
(\bF^*_{q^l}/\bF^{*}_q)\times(\bF^*_{q}/S_{\bF_q})$. Let $\eta_{\bF^*_{q}}$ be the non-principal character of $\bF^*_{q}/S_{\bF_q}$. Then $\eta_{\bF^*_{q}}$ can be extended to a
ring homomorphism $\bZ[\bF^*_{q^l}/S_{\bF_q}]\to\bZ[\bF^*_{q^l}/\bF^{*}_q]$, which will be again denoted by $\eta_{\bF^*_{q}}$. The ring homomorphism
$\eta_{\bF^*_{q}}$ amounts to replace the non-identity element of $\bF^*_{q}/S_{\bF_q}$ by $-1$ and clearly $\eta_{\bF^*_{q}}(R)=W_{q^l/q}$. Also $\eta_{\bF^*_{q}}$ induces an one-to-one
correspondence between transversals of $\bF^*_{q}/S_{\bF_q}$ in $\bF^*_{q^l}/S_{\bF_q}$ and elements in $\bZ[\bF^*_{q^l}/\bF^{*}_q]$ with $\pm 1$
coefficients, which we denote by $\tilde\eta_{\bF^*_{q}}$.

{\cor Let $l$ be odd, $D$ be an element with $\pm 1$ coefficients in
$\bZ[\bF^*_{q^l}/\bF^{*}_q]$ and $W_{q^l/q}$ be the
Singer circulant weighing matrix in $\bZ[\bF^*_{q^l}/\bF^{*}_q]$.
Then the following statements are equivalent:
\begin{itemize}
\item[(i)] $(\tilde\eta_{\bF^*_{q}}\circ\mu)^{-1}(D)$ is a Paley type group scheme
in $\bF_{q^l}$;
\item[(ii)] there is another
$\pm 1$ coefficient element $\widehat{D}$ in
$\bZ[\bF^*_{q^l}/\bF^{*}_q]$ such that
$$D^{(-1)}W_{q^l/q}=q^{(l-1)/2}\widehat{D}$$ in
$\bZ[\bF^*_{q^l}/\bF^{*}_q]$;
\item[(iii)] $D^{(-1)}W_{q^l/q}$ is divisible
by $q^{(l-1)/2}$ in $\bZ[\bF^*_{q^l}/\bF^{*}_q]$.\end{itemize}}\vspace{2mm}

For example, when $D=\bF^*_{q^l}/\bF^{*}_q$, one obtains the classical Paley group scheme $S_{\bF_{q^l}}$
 in $\bF_{q^l}$. In this case, the element
$\widehat{D}$ is either $D$ or $-D$ depending on the value of the sum of
coefficients of the Singer weighing matrix $W_{q^l/q}$. For each element $D$ with $\pm 1$ coefficients in $\bZ[\bF^*_{q^l}/\bF^{*}_q]$, there are
two subsets $D_+$ and $D_-$ in $\bF^*_{q^l}/\bF^{*}_q$ such that $D_++D_-=\bF^*_{q^l}/\bF^{*}_q$, $D_+-D_-=D$, and by equation~(\ref{DX}), $(\tilde\eta_{\bF^*_{q}}\circ\mu)^{-1}(D)=\bD(D_+)$.
Since $D^{(-1)}W_{q^l/q}$ is divisible by $q^{(l-1)/2}$ in $\bZ[\bF^*_{q^l}/\bF^{*}_q]$ if and only if $D_+^{(-1)}W_{q^l/q}$ is divisible by $q^{(l-1)/2}$ in $\bZ[\bF^*_{q^l}/\bF^{*}_q]$,
and $D$ is a projective half-point set in $\bF^*_{q^l}$ over $\bF^{*}_q$ if and only if $D=\bD(X)$ for some subset $X$ in $\bF^*_{q^l}/\bF^{*}_q$ when $l$ is odd, we have

{\thm\label{div} Let $l$ be odd and $X$ be a subset of $\bF^*_{q^l}/\bF^{*}_q$. Then $\bD(X)$ is a Paley type group scheme in $\bF_{q^l}$ if and only if $X^{(-1)}W_{q^l/q}$ is divisible by $q^{(l-1)/2}$ in $\bZ[\bF^*_{q^l}/\bF^{*}_q]$.}

{\rem For any odd integer $l$ and any $S_{\bF_q}$ invariant Paley type group scheme $\bD$ in $\bF_{q^l}$, the set $X$ is actually given by $X=\gamma(\bD\cap S_{\bF_{q^l}})$ and $\bD=\bD(X)$, where $\gamma:\bF^*_{q^l}\to\bF^*_{q^l}/\bF^*_q$ is the natural projection. It measures how much $\bD$ remains the same as or how much $\bD$ deviates from the standard Paley group schemes $S_{\bF_{q^l}}$ and $N_{\bF_{q^l}}$ in $\bF_{q^l}$.}

\section{Proofs of Theorems~\ref{main2}--\ref{main3}}\label{proof}
In order to prove Theorems~\ref{main2}--\ref{main3}, we need to use Gauss sums over finite fields because they are related to the character sums of Singer difference sets and Singer weighing matrices. Let $\xi_p$ be a primitive $p$-th root of unity in the complex number field $\bC$. For each character $\chi\in\widehat{\bF_q^*}$, the Gauss sum $G_{\bF_q}(\chi)$ of $\chi$ over $\bF_q$ is defined to be
$$G_{\bF_q}(\chi)=\sum_{x\in\bF_q^*}\chi(x)\xi^{\tr_{q/p}(x)}.$$

In \cite{Y}, Yamamoto proved the following lemma.

{\lem\cite{Y}\label{moto} For the Singer difference set $S_{q^l/q}$ in $\bF_{q^l}^*/\bF_q^*$, we have $$\chi(S_{q^l/q})=-\frac{G_{\bF_{q^l}}(\chi)}{q}$$ for each non-principal character $\chi\in\widehat{\bF_{q^l}^*/\bF_q^*}$.}

Let $\chi\in\widehat{\bF^*_{q^l}/\bF^*_q}$ and $\eta_{\bF^*_{q^l}}$ be the quadratic character of $\bF^*_{q^l}$. Since $l$ is odd, the restriction $\eta_{\bF^*_{q^l}}|_{\bF^*_q}=\eta_{\bF^*_{q}}$. Therefore the Gauss sum

\begin{align*}G_{\bF_{q^l}}(\chi\eta_{\bF^*_{q^l}})&=\sum_{x\in\bF^*_{q^l}}\xi_p^{\tr_{q^l/p}(x)}(\chi\eta_{\bF^*_{q^l}})(x)\\
                         &=\sum_{x\in\bF^*_{q^l}}\xi_p^{\tr_{q/p}(\tr_{q^l/q}(x))}\chi(x)\eta_{\bF^*_{q^l}}(x)\\
                         &=\sum_{x\in H^*}\xi_p^{\tr_{q/p}(\tr_{q^l/q}(x))}\chi(x)\eta_{\bF^*_{q^l}}(x)+\sum_{x\in\bF^*_q}\sum_{y\in S_{q^l/q}}\xi_p^{\tr_{q/p}(\tr_{q^l/q}(xy))}\chi(xy)\eta_{\bF^*_{q^l}}(xy)\\
                         &=\sum_{x\in H^*}\chi(x)\eta_{\bF^*_{q^l}}(x)+\sum_{x\in\bF^*_q}\sum_{y\in S_{q^l/q}}\xi_p^{\tr_{q/p}(x)}\eta_{\bF^*_{q}}(x)\chi(y)\eta_{\bF^*_{q^l}}(y)\\
                         &=G_{\bF_q}(\eta_{\bF^*_{q}})(\chi\eta_{\bF^*_{q^l}})(S_{q^l/q})\\
                         &=G_{\bF_q}(\eta_{\bF^*_{q}})\chi(W_{q^l/q}),\\
\end{align*} where $H^*=\{x\in\bF_{q^l}^*~|\tr_{q^l/q}(x)=0\}$ as in the proof of Theorem~\ref{mul1}.
By the Davenport-Hasse product formula (see~\cite{BEW}) and Lemma~\ref{moto}, we have
$$G_{\bF_{q^l}}(\chi\eta_{\bF^*_{q^l}})=\frac{G_{\bF_{q^l}}(\chi^2)G_{\bF_{q^l}}(\eta_{\bF^*_{q^l}})}{\chi(2)^2G_{\bF_{q^l}}(\chi)}$$
and $$\chi(W_{q^l/q})=\frac{G_{\bF_{q^l}}(\chi\eta_{\bF^*_{q^l}})}{G_{\bF_q}(\eta_{\bF^*_{q}})}=\chi(4^{-1})\frac{G_{\bF_{q^l}}(\eta_{\bF^*_{q^l}})}{G_{\bF_q}(\eta_{\bF^*_{q}})}\frac{G_{\bF_{q^l}}(\chi^2)}{G_{\bF_{q^l}}(\chi)}
=\pm\chi(4^{-1})q^{\frac{l-1}{2}}\frac{\chi(S_{q^l/q}^{(2)})}{\chi(S_{q^l/q})}.$$

Hence Theorem~\ref{div} can now be restated by using characters of $\bF^*_{q^l}/\bF^*_q$.

{\thm\label{cdiv} Let $l$ be odd and $X$ be a subset of $\bF^*_{q^l}/\bF^{*}_q$. Then $\bD(X)$ is a Paley type group scheme in $\bF_{q^l}$ if and only if $\chi(S_{q^l/q}^{(2)})$ divides $\chi(X)\chi(S_{q^l/q})$ for all non-principal $\chi\in\widehat{\bF_{q^l}^*/\bF_q^*}$.}
\proof By Theorem~\ref{div}, $\bD(X)$ is a Paley type group scheme in $\bF_{q^l}$ if and only if $X^{(-1)}W_{q^l/q}$ is divisible by $q^{(l-1)/2}$ in $\bZ[\bF^*_{q^l}/\bF^{*}_q]$.  By the Fourier inversion formula, $X^{(-1)}W_{q^l/q}$ is divisible by $q^{(l-1)/2}$ in $\bZ[\bF^*_{q^l}/\bF^{*}_q]$ if and only if $q^{(l-1)/2}$ divides $\chi(X^{(-1)})\chi(W_{q^l/q})$ for every $\chi\in\widehat{\bF_{q^l}^*/\bF_q^*}$, as $q^{(l-1)/2}$ is prime to $(q^t-1)/(q-1)$. From $\chi(W_{q^l/q})=\pm\chi(4^{-1})q^{\frac{l-1}{2}}{\chi(S_{q^l/q}^{(2)})}/{\chi(S_{q^l/q})}$ and $\chi(S_{q^l/q})\chi(S_{q^l/q}^{(-1)})=\chi(S_{q^l/q}^{(2)})\chi(S_{q^l/q}^{(-2)})$, we find that
$q^{(l-1)/2}$ divides $\chi(X^{(-1)})\chi(W_{q^l/q})$ if and only if $\chi(X)\chi(S_{q^l/q})$ is divisible by $\chi(S_{q^l/q}^{(2)})$ for all $\chi\in\widehat{\bF_{q^l}^*/\bF_q^*}$.\qed

As an application of Theorem~\ref{cdiv}, we give a new proof of Theorem 3.6 and Corollary 3.7 in~\cite{FMX}, which slightly generalizes these results. These results themselves are generalizations of Theorem 3.6 in~\cite{FX}.
To this end, we need the following proposition of Langevin~\cite{La}.

{\prop\label{la}~\cite[Proposition 4.2]{La} Let $m$ be a positive integer and $p'\ne 3$ be a prime such that $p'\equiv 3\pmod 8$ and the order $l$ of $p$ in $\bZ^*_{{p'}^m}$ is $l=({p'}^m-{p'}^{m-1})/2$. Let $q=p^l$ and $\chi\in\widehat{\bF_q^*/\bF_p^*}$ such that the order of $\chi$ is ${p'}^m$. Then $$G_{\bF_q}(\chi)=p^{\frac{l-h}{2}}\frac{a+b\sqrt{-p'}}{2},$$ where $h$ is the class number of $\bQ(\sqrt{-p'})$ and $a$ and $b$ are integers such that \begin{itemize}
\item[(i)] $p$ does not divide $b$,
\item[(ii)] $a\equiv -2p^{(l+h)/2}\pmod {p'}$,
\item[(iii)] $a^2+b^2p'=4p^h$.\end{itemize}}
Items (i) and (ii) imply that $ab\ne0$. Therefore when $4p^h=1+p'$, the Gauss sum $$G_{\bF_q}(\chi)=p^{\frac{l-h}{2}}\frac{\pm1\pm\sqrt{-p'}}{2}.$$ Since $\frac{\pm1\pm\sqrt{-p'}}{2}\in\{\chi'(S_{\bF_{p'}}),\chi'(N_{\bF_{p'}}),\chi'(S_{\bF_{p'}}\cup\{0\}),\chi'(N_{\bF_{p'}}\cup\{0\})\}$ for every non-principal character $\chi'$ of the additive group of $\bF_{p'}$, and $p$ generates $S_{\bF_{p'}}$ as the order of $p$ in $\bZ^*_{{p'}^m}$ is $({p'}^m-{p'}^{m-1})/2$, Lemma~\ref{moto} and Proposition~\ref{la} imply the following corollary.

{\cor\label{corla} Let $m$ be a positive integer and $p'\ne 3$ be a prime such that $p'\equiv 3\pmod 8$, the order $l$ of $p$ in $\bZ^*_{{p'}^m}$ is $l=({p'}^m-{p'}^{m-1})/2$ and $4p^h=1+p'$,  where $h$ is the class number of $\bQ(\sqrt{-p'})$. Let $q=p^l$. If we identify the subgroup of order $p'$ in $\bZ_{{p'}^m}$ with $\bF_{p'}$, then there is a unique subset $P_{q/p}\in\{S_{\bF_{p'}},N_{\bF_{p'}},S_{\bF_{p'}}\cup\{0\},N_{\bF_{p'}}\cup\{0\}\}$ in $\bF_{p'}$ such that for every $\chi\in\widehat{\bF_q^*/\bF_p^*}$ of order ${p'}^m$, $$\chi(S_{q/p})=-p^{\frac{l-h-2}{2}}\chi(P_{q/p}).$$}

Corollary~\ref{corla} yields the following theorem which is slightly more general than Theorem 3.6 and Corollary 3.7 in~\cite{FMX} because our $T$ in the theorem is an arbitrary transversal.

{\thm\label{slight} Let $m$ be a positive integer and $p'\ne 3$ be a prime such that $p'\equiv 3\pmod 8$, the order $l$ of $p$ in $\bZ^*_{{p'}^m}$ is $l=({p'}^m-{p'}^{m-1})/2$ and $4p^h=1+p'$,  where $h$ is the class number of $\bQ(\sqrt{-p'})$. Let $q=p^l$ and $\gamma:\bF_q/\bF_p\to\bZ_{{p'}^m}$ be the natural projection. Let $P_{q/p}$ be the subset in the subgroup of order $p'$ in $\bZ_{{p'}^m}$ as in Corollary~\ref{corla}. Then for every transversal $T$ of the subgroup of order $p'$ in $\bZ_{{p'}^m}$, $TP_{q/p}^{(2)}$ is a subset in $\bZ_{{p'}^m}$ and $\bD(\gamma^{-1}(TP_{q/p}^{(2)}))$ is a Paley type group scheme in $\bF_q$. }
\proof Let $\chi$ be a non-principal character of $\bF^*_{q}/\bF^*_p$. If $\chi$ is non-principal on the kernel of $\gamma$,
 then $\chi(\gamma^{-1}(TP_{q/p}^{(2)}))\chi(S_{q/p})=0$ and $\chi(S_{q/p}^{(2)})$ divides $\chi(\gamma^{-1}(TP_{q/p}^{(2)}))\chi(S_{q/p})$. If $\chi$ is principal on the kernel of $\gamma$, then $\chi\in\widehat{\bZ}_{{p'}^m}$. If $\chi$ has order ${p'}^m$, then by Corollary~\ref{corla}, we have that the character sum $\chi(S_{q/p}^{(2)})=\chi^2(S_{q/p})=-p^{\frac{l-h-2}{2}}\chi^2(P_{q/p})=-p^{\frac{l-h-2}{2}}\chi(P_{q/p}^{(2)})$ divides $$\chi(\gamma^{-1}(TP_{q/p}^{(2)}))\chi(S_{q/p})=-p^{\frac{l-h-2}{2}}|\Ker(\gamma)|\chi(T)\chi(P^{(2)}_{q/p})\chi(P_{q/p}).$$ If $\chi$ has order dividing ${p'}^m$ but not equal to ${p'}^m$, then $\chi(T)=0$ as $T$ is a transversal of the subgroup of order $p'$ in $\bZ_{{p'}^m}$. Hence $\chi(S_{q/p}^{(2)})$ divides $\chi(\gamma^{-1}(TP_{q/p}^{(2)}))\chi(S_{q/p})=0$. By Theorem~\ref{cdiv}, $\bD(\gamma^{-1}(TP_{q/p}))$ is a Paley type group scheme in $\bF_q$.\qed

We now prove Theorems~\ref{main2}--\ref{main3}.

{\it Proof of Theorem~\ref{main2}}: Let $X$ be a $((q^l-1)/(q-1),q^{l-1},q^{l-2}(q-1))$-difference set in $\bF^*_{q^l}/\bF^*_q$. If $\bD(X)$ is a Paley type group scheme in $\bF_{q^l}$, then by Theorem~\ref{cdiv},
$\chi(S_{q^l/q}^{(2)})$ divides $\chi(X)\chi(S_{q^l/q})$ for every $\chi\in\widehat{\bF_{q^l}^*/\bF_q^*}$, and therefore $q^{l-2}$ divides $\chi(X)\chi(S_{q^l/q})\chi(S_{q^l/q}^{(-2)})$, or equivalently, $q^{l-2}$ divides $\chi(X^{(-1)})\chi(S^{(-1)}_{q^l/q})\chi(S_{q^l/q}^{(2)})$ for every $\chi\in\widehat{\bF_{q^l}^*/\bF_q^*}$. By the Fourier inversion formula and Theorem~\ref{product}, there is a $((q^l-1)/(q-1),q^{l-1},q^{l-2}(q-1))$-difference set $Y$ in $\bF^*_{q^l}/\bF^*_q$ such that $$X^{(-1)}S^{(-1)}_{q^l/q}S_{q^l/q}^{(2)}=(q^{l-2}+(q^{l-1}-q^{l-2})\bF^*_{q^l}/\bF^*_q)Y.$$ This implies that $XY=S^{(-1)}_{q^l/q}S_{q^l/q}^{(2)}$ and $X$ is an Arasu-Dillon-Player difference set.

Conversely, if $X$ is an Arasu-Dillon-Player difference set, then by Remark~\ref{remadp}, $S^{(2)}_{q^l/q}$ divides $XS_{q^l/q}$, and therefore $\chi(S_{q^l/q}^{(2)})$ divides $\chi(X)\chi(S_{q^l/q})$ for every $\chi\in\widehat{\bF_{q^l}^*/\bF_q^*}$. By Theorem~\ref{cdiv}, $\bD(X)$ is a Paley type group scheme in $\bF_{q^l}$.\qed

{\it Proof of Theorem~\ref{main1}}: Let $\chi$ be a non-principal character of $\bF^*_{q^l}/\bF^*_l$. If $\chi$ is non-principal on the kernel of $\gamma$,
 then $\chi(\gamma^{-1}(X))\chi(S_{q^l/q})=0$ and $\chi(S_{q^l/q}^{(2)})$ divides $\chi(\gamma^{-1}(X))\chi(S_{q^l/q})$. If $\chi$ is principal on the kernel of $\gamma$,  then $\chi\in\widehat{\bZ}_n$ and $\chi(S_{q^l/q}^{(2)})=\chi^2(S_{q^l/q})=\chi^{p^t}(S_{q^l/q})=\chi(S_{q^l/q}^{(p^t)})$ for some integer $t$ as $2\in\<p\>$ in $\bZ_n^*$, and $\chi(S_{q^l/q}^{(2)})=\chi(S_{q^l/q})$ as $p\in\cM_0(S_{q^l/q})$. This again implies that $\chi(S_{q^l/q}^{(2)})$ divides $\chi(\gamma^{-1}(X))\chi(S_{q^l/q})$. By Theorem~\ref{cdiv}, $\bD(\gamma^{-1}(X))$ is a Paley type group scheme in $\bF_{q^l}$. \qed

{\it Proof of Theorem~\ref{main3}}: If $\bD(X)$ is a Paley type group scheme in $\bF_{q^t}$, by Theorem~\ref{cdiv}, $\chi(S^{(2)}_{q^t/q})$ divides $\chi(X)\chi(S_{q^t/q})$ for all $\chi\in\widehat{\bF_{q^t}^*/\bF_q^*}$. By Gordon-Mills-Welch decomposition, $S^{(2)}_{q^{st}/q}=\tilde{R}^{(2)}_{q^{st}/q^t}S^{(2)}_{q^{t}/q}$ and $S_{q^{st}/q}=\tilde{R}_{q^{st}/q^t}S_{q^{t}/q}$. Therefore $\chi(S^{(2)}_{q^{st}/q})=\chi(\tilde{R}^{(2)}_{q^{st}/q^t})\chi(S^{(2)}_{q^{t}/q})$ divides $\chi(\tilde{R}^{(2)}_{q^{st}/q^t})\chi(X)\chi(S_{q^t/q})$, which also divides $\chi(\tilde{R}^{(2)}_{q^{st}/q^t})\chi(X)\chi(\tilde{R}_{q^{st}/q^t})\chi(S_{q^t/q})=\chi(\tilde{R}^{(2)}_{q^{st}/q^t}X)\chi(S_{q^{st}/q^t})$ for all $\chi\in\widehat{\bF_{q^{st}}^*/\bF_q^*}$. By Theorem~\ref{cdiv},
$\bD(\tilde{R}^{(2)}_{q^{st}/q^t}X)$ is a Paley type group scheme in $\bF_{q^{st}}$. Since $\tilde{R}^{(2)}_{q^{st}/q^t}\tilde{R}^{(-2)}_{q^{st}/q^t}=\tilde{R}_{q^{st}/q^t}\tilde{R}^{(-1)}_{q^{st}/q^t}$,
$\chi(\tilde{R}^{(2)}_{q^{st}/q^t})$ divides $\chi(\tilde{R}_{q^{st}/q^t})\chi(\tilde{R}^{(-1)}_{q^{st}/q^t})$ for all $\chi\in\widehat{\bF_{q^{st}}^*/\bF_q^*}$. Therefore $\chi(S^{(2)}_{q^{st}/q})=\chi(\tilde{R}^{(2)}_{q^{st}/q^t})\chi(S^{(2)}_{q^{t}/q})$ divides $\chi(\tilde{R}_{q^{st}/q^t})\chi(\tilde{R}^{(-1)}_{q^{st}/q^t})\chi(X)\chi(S_{q^{t}/q})=\chi(\tilde{R}^{(-1)}_{q^{st}/q^t}X)\chi(S_{q^{st}/q})$ for all $\chi\in\widehat{\bF_{q^{st}}^*/\bF_q^*}$, and by Theorem~\ref{cdiv},
$\bD(\tilde{R}^{(-1)}_{q^{st}/q^t}X)$ is a Paley type group scheme in $\bF_{q^{st}}$. \qed

\section{Conclusions}\label{con}
In $\bF_{5^3}$, the Paley type group scheme $\bD(S^{(2)}_{5^3/5})$ has $|\Aut(\mathfrak{C}(\bD(S^{(2)}_{5^3/5}))|=2^3\cdot3\cdot5^3$ and Theorem~\ref{main2} replaces one of the question marks in Table~\ref{t1}.

Using Theorem~\ref{adp} and Theorem~\ref{main2}, we did a MAGMA search and found the following 10 inequivalent Paley type group schemes in $\bF_{3^5}$ which have non-isomorphic configurations:
$\bD(S^{(2)}_{3^5/3})$, $\bD(S^{(4)}_{3^5/3})$, $\bD(S^{(5)}_{3^5/3})$, $\bD(S^{(10)}_{3^5/3})$, $\bD(S^{(20)}_{3^5/3})$, $\bD(S^{(40)}_{3^5/3})$, $\bD(A_{3^5/3}(4))$, $\bD(A_{3^5/3}(5))$,
$\bD(A_{3^5/3}(10))$, $\bD(A_{3^5/3}(20))$. These Paley type group schemes are all $\Gal(\bF_{3^5})$ invariant and none of them is equivalent to the Paley type group scheme from the 3 semifields mentioned in~\cite{CK}. The scheme $\bD(S_{3^5/3}^{(10)})$ is equivalent to ${\rm RT}(1)$ while $\bD(A_{3^5/3}(10))$ is equivalent to ${\rm RT}(-1)$ (see~\cite{DW}). Theorem~\ref{div} or Theorem~\ref{cdiv} clearly has the following consequence.
{\thm\label{dissum} Let $l$ be an odd integer and $X_1$ and $X_2$ be subsets of $\bF^*_{q^l}/\bF^{*}_q$. If $X_1\cap X_2=\emptyset$ and $\bD(X_1)$ and $\bD(X_2)$ are both Paley type group schemes in $\bF_{q^l}$, then $\bD(X_1\cup X_2)$ is also a Paley type group scheme in $\bF_{q^l}$. Equivalently, If $X_1\cup X_2=\bF^*_{q^l}/\bF^{*}_q$ and $\bD(X_1)$ and $\bD(X_2)$ are both Paley type group schemes in $\bF_{q^l}$, then $\bD(X_1\cap X_2)$ is also a Paley type group scheme in $\bF_{q^l}$.}

Among the 10 Arasu-Dillon-Player difference sets we used in $\bF_{3^5}^*/\bF_3^*$, we found that $S^{(5)}_{3^5/3}\cup A_{3^5/3}(5)=\bF_{3^5}^*/\bF_3^*$ and, by Theorem~\ref{dissum}, $\bD(S^{(5)}_{3^5/3}\cap A_{3^5/3}(5))$ is a
Paley type group scheme, whose configuration is not isomorphic to any of the configurations of the aforementioned 13 Paley type group schemes in $\bF_{3^5}$.
We now understand how to construct 14 of the 58 $\Gal(\bF_{3^5})$ invariant Paley type group schemes in Table~\ref{t2} and still have 44 more to go.

In $\bF_{7^3}$, Theorem~\ref{main2} yields two non-isomorphic configurations from $\bD(S^{(2)}_{7^3/7})$ and $\bD(S^{(-1)}_{7^3/7})$ and these two Paley type group schemes have
$$|\Aut(\mathfrak{C}(\bD(S^{(2)}_{7^3/7}))|=|\Aut(\mathfrak{C}(\bD(S^{(-1)}_{7^3/7}))|=3^2\cdot7^3.$$

In $\bF_{3^7}$, we used a MAGMA program and found that $\bD(S^{(2)}_{3^7/3})$, $\bD(S^{(4)}_{3^7/3})$, $\bD(S^{(5)}_{3^7/3})$, $\bD(S^{(10)}_{3^7/3})$, $\bD(S^{(14)}_{3^7/3})$, $\bD(S^{(28)}_{3^7/3})$, $\bD(S^{(182)}_{3^7/3})$, $\bD(S^{(364)}_{3^7/3})$, $\bD(A_{3^7/3}(4))$, $\bD(A_{3^7/3}(5))$, $\bD(A_{3^7/3}(10))$, $\bD(A_{3^7/3}(14))$, $\bD(A_{3^7/3}(28))$ and $\bD(A_{3^7/3}(182))$ are all the inequivalent Paley type group schemes with non-isomorphic configurations that can be obtained from Theorem~\ref{main2} and Theorem~\ref{adp}. None of the configurations of these Paley type group schemes is isomorphic to that of ${\rm DY}(1)$, ${\rm DY}(-1)$, ${\rm RT}(1)$ or ${\rm RT}(-1)$, where ${\rm DY}(\pm1)$ are the Paley type group schemes
constructed in ~\cite{DY} and ${\rm RT}(\pm1)$ are those in~\cite{DW}. Therefore, besides Paley group scheme, there are at least 18 $\Gal(\bF_{3^7})$ invariant Paley type group schemes with non-isomorphic configurations in $\bF_{3^7}$.

\vspace{3mm}

{\bf Acknowledgement:} Y. Q. Chen would like to thank the Department of Mathematics at Zhejiang University for the hospitality he received during his visit when this research was initiated.
The work of T. Feng was supported in part by the Fundamental Research Funds for the Central Universities, Zhejiang Provincial Natural Science Foundation.


\end{document}